\documentclass[12pt,a4paper]{article}
\usepackage{amssymb}

\usepackage{graphicx,amssymb,amsfonts,epsfig,amsthm,a4,amsmath,url}
\usepackage[latin1]{inputenc}

\newtheorem{Thm}{Theorem}[section]
\newtheorem{Prop}[Thm]{Proposition}
\newtheorem{Lem}[Thm]{Lemma}
\newtheorem{Cor}[Thm]{Corollary}
\newtheorem{Ex}[Thm]{Example}
\theoremstyle{definition}
\newtheorem{Def}[Thm]{Definition}

\newtheorem{Rem}[Thm]{Remark}

\newcommand{\Z}{\mathbf{Z}}

\newcommand{\N}{\mathbf{N}}
\newcommand{\R}{\mathbf{R}}

\newcommand{\Q}{\mathbf{Q}}

\newcommand{\T}{\mathbf{T}}

\newcommand{\HH}{\mathcal{H}}

\title{Locally compact groups with every isometric action bounded or proper}
\author{Romain TESSERA and Alain VALETTE\\ {\small (with an appendix by Corina CIOBOTARU)}}
\date{\today}

\begin{document}

\baselineskip=16pt

\maketitle

\begin{abstract} A locally compact group $G$ has property PL if every isometric $G$-action either has bounded orbits or is (metrically) proper. For $p>1$, say that $G$ has property $BP_{L^p}$ if the same alternative holds for the smaller class of affine isometric actions on $L^p$-spaces. We explore properties PL and $BP_{L^p}$ and prove that they are equivalent for some interesting classes of groups: abelian groups, amenable almost connected Lie groups, amenable linear algebraic groups over a local field of characteristic 0.

The appendix provides new examples of groups with property PL, including non-linear ones. 


\end{abstract}

\section{Introduction}

Let the locally compact group $G$ act by isometries on a metric space $(X,d)$. The action is {\it locally bounded} if $Kx$ is bounded for every $x\in X$ and every compact set $K\subset G$; the action is {\it bounded} if every orbit is bounded. On the other hand, the action is (metrically) {\it proper} if $\lim_{g\rightarrow\infty} d(gx,x)=+\infty$ for every $x\in X$.

A {\it length function} on $G$ is a non-negative function $L$ on $G$ which is bounded on compact subsets, is symmetric ($L(g)=L(g^{-1})$ for every $g\in G$), and is sub-additive: $L(gh)\leq L(g)+L(h)$ for every $g,h\in G$. Clearly if $G$ admits a locally bounded action by isometries on $(X,d)$, then for every $x\in X$ the function $L:G\rightarrow\R^+:g\mapsto d(gx,x)$ is a length function on $G$. It is known that $G$ admits a proper length function if and only if $G$ is $\sigma$-compact (see Section 2 in \cite{Cor}). In the next definition, the equivalence of the two statements is Proposition 1.2 in \cite{Cor}.

\begin{Def} (see \cite{Cor}) A locally compact group $G$ has property PL if every locally bounded action of $G$ by isometries is either bounded or proper; equivalently, every length function on $G$ is either bounded or proper.
\end{Def}

For $p\geq 1$, a length function $L$ on $G$ is {\it a $L^p$-type length function} if it comes from a continuous affine isometric action $\alpha$ of $G$ on some $L^p$-space $L^p(X,\mu)$, i.e. $L(g)=\|\alpha(g)x-x\|$ for some $x\in L^p(X,\mu)$. In the terminology of Definition 6.5 in \cite{CDH}, the invariant kernel $(g,h)\mapsto L(g^{-1}h)^p$ has type $p$ on $G$.

\begin{Def}
For $p\geq 1$, a locally compact group $G$ has property $BP_{L^p}$ if every affine isometric action of $G$ on a $L^p$-space is either bounded or proper; equivalently every $L^p$-type length function on $G$ is either bounded or proper.
\end{Def}

Recall from \cite{BFGM} that $G$ has property $FL^p$ if every continuous affine isometric action of $G$ on a $L^p$-space, has a fixed point\footnote{Property $FL^2$ is more commonly denoted by FH and, for $\sigma$-compact locally compact groups, property FH is equivalent to Kazhdan's property (T); see \cite{BHV} for all this.}. Obviously property $BP_{L^p}$ is implied both by property PL and by property $FL^p$. 
 
A surprising fact, discovered by Y. Shalom (\cite{ShaAnn}, Theorem 3.4), is that simple Lie groups with finite center have property $BP_{L^2}$. Since those have property $FH$ except when locally isomorphic to $\mathrm{SO}(n,1)$ or $\mathrm{SU}(n,1)$, this is really a statement about the latter two classes of groups. A stronger result was proved by Y. Cornulier (\cite{Cor}, Theorem 1.4): property PL holds for all simple algebraic groups over a local field\footnote{For isometric actions which are {\it continuous}, not just locally bounded, an even stronger result holds for simple algebraic groups over local field: a continuous isometric action either is proper or has a globally fixed point, see Theorem 6.1 in \cite{BadGel}.}. In the same paper, it is also proved that certain semi-direct products have property PL, e.g. $\R^d\rtimes K$, where $K$ is a closed subgroup of the orthogonal group $\mathrm{O}(d)$ acting transitively on the unit sphere (see Proposition 1.8 in \cite{Cor}); or $F\rtimes A^\times$, where $F$ is a non-Archimedean local field and $A^\times$ is the invertible group of its ring of integers (Proposition 1.9 in \cite{Cor}).

The aim of the present paper is to investigate the relation between properties PL and $BP_{L^p}$. It was a surprise for us that, for some interesting classes of groups, they turn out to be equivalent. For example, for abelian groups, both are equivalent to compactness:

\begin{Thm}\label{abelian} Let $A$ be a locally compact abelian (LCA) group. The following are equivalent:
\begin{enumerate} 
\item[a)] $A$ has property PL;
\item[b)] for some (resp.\ every) $p\geq1$ the group $A$ has property $BP_{L^p}$;
\item[c)] $A$ is compact.
\end{enumerate}
\end{Thm}

For amenable locally compact groups, we have:

\begin{Thm}\label{semidir} Let $G$ be a locally compact group admitting a closed co-compact normal subgroup $V$ such that $V\simeq F^d$ with $F$ a local field of characteristic 0 and $d>0$ (so that the compact group $G/V$ acts on $V$). The following are equivalent:
\begin{enumerate}
\item[a)] $G$ has property PL;
\item[b)] for some (resp.\ every) $p\geq 1$ the group $G$ has property $BP_{L^p}$;
\item[c)] $G/V$ is infinite and it acts irreducibly on $V$.
\end{enumerate}
\end{Thm}

{\bf Notation:} We denote by $\mathcal{G}_{LA}$ the union of the class of almost connected Lie groups and the class of groups of the form $\mathbb{G}(F)$, the group of $F$-rational points of a linear algebraic group $\mathbb{G}$ defined over a non-Archimedean local field $F$ with characteristic 0.

\begin{Thm}\label{algebraicamenable} Let $G$ be an amenable non-compact group in $\mathcal{G}_{LA}$. The following are equivalent:
\begin{enumerate}
\item[a)] $G$ has property PL;
\item[b)] for some (resp.\ every) $p\geq 1$ the group $G$ has property $BP_{L^p}$;
\item[c)] there exists a compact normal subgroup $W$ of $G$ such that $H:=G/W$ has a closed co-compact subgroup $V$ isomorphic to $F^d$ for some $d\geq 1$, with $H/V$ infinite and acting irreducibly on $V$.
\end{enumerate}
\end{Thm}


For a group $G$, we denote by $Ad(G)$ the image of $G$ in its group of inner automorphisms. For non-amenable groups we have:

\begin{Thm}\label{algebraicPL} Let $G$ be a non-amenable locally compact group which is either an almost connected Lie group or a linear algebraic group over a local field of any characteristic. The following are equivalent:
\begin{enumerate}
\item[a)] $G$ has property PL;
\item[b)] every closed normal subgroup of $G$ is either compact or co-compact;
\item[c)] there exists a compact normal subgroup $W$ of $G$ such that $G/W$ admits a closed, co-compact, normal subgroup $N$ which is isomorphic to a direct product $S_1\times...\times S_n$ of simple groups, and the simple factors $S_1,...,S_n$ are permuted transitively under $Ad(G)$.
\end{enumerate}
\end{Thm}
The previous result actually holds under weaker assumptions on $G$, see Theorem \ref{P.-E.C.} for the precise statement.

The linear algebraic groups with property FH have been characterized by S.P. Wang \cite{Wang}. So to classify linear algebraic groups with property $BP_{L^2}$, we may assume that they do not have property FH.

\begin{Thm}\label{algebraic} Let $G$ be a group in $\mathcal{G}_{LA}$. Assume that $G$ does NOT have property FH and is non-amenable. The following are equivalent:
\begin{enumerate}
\item[a)] The group $G$ has property $BP_{L^2}$.
\item[b)] $G$ admits a finite normal subgroup $W$ such that $G/W$ admits a closed, co-compact, normal subgroup $N$ which is isomorphic to a direct product $S_1\times...\times S_n$ of simple groups, and the simple factors $S_1,...,S_n$ are permuted transitively under $Ad(G)$.
Moreover, if $G=\mathbb{G}(F)$ with $F$ is non-Archimedean, each simple factor of $N$ is a simple algebraic group of rank 1 over $F$; if $G$ is Lie almost connected, each simple factor of $N$ is locally isomorphic to $\mathrm{SO}(n,1)$ or $\mathrm{SU}(m,1)$.
\end{enumerate}
\end{Thm}

Finally, we prove a general result about centers of $BP_{L^p}$-groups.

\begin{Thm}\label{centralext} Fix $p>1$. Let $G$ be a compactly generated locally compact group satisfying property $BP_{L^p}$ but not property $FL^p$. Then the center of $G$ is compact.
\end{Thm}

In a previous paper \cite{CTV}, property $BP_0$ was introduced for a locally compact group $G$: it means that $G$ satisfies the bounded/proper alternative for affine isometric actions on Hilbert spaces, such that the linear part is a $C_0$, or mixing, representation. The class of groups with $BP_0$ is significantly larger than the class of groups with $BP_{L^2}$. For example, it was proved in \cite{CTV} that every group with non-compact center (in particular every abelian group) has $BP_0$.

The paper is structured as follows. Section 2 contains generalities on property $BP_{L^p}$. In particular we prove that, for $1\leq p\leq 2$ and $G$ locally compact separable, property $BP_{L^p}$ is equivalent to every action of $G$ on a space with measured walls being bounded or proper (Proposition \ref{measuredwalls}). Section 3 contains generalities on property PL. Theorems \ref{abelian}, \ref{semidir}, \ref{algebraicamenable} and \ref{centralext} are proved in section 4, which is the core of the paper. Theorem \ref{algebraicPL} is proved in section 5. Section 6 deals specifically with property $BP_{L^2}$: we prove Theorem \ref{algebraic} and make in Proposition \ref{PropBek} the connection with the Howe-Moore property, by proving that it implies property $BP_{L^2}$. This provides a new proof of the already mentioned Theorem 3.4 in \cite{ShaAnn}, stating that $\mathrm{SO}(n,1)$ and $\mathrm{SU}(n,1)$ have property $BP_{L^2}$ (the original proof used the Mautner phenomenon). Since property $BP_{L^2}$ is implied both by property PL and the Howe-Moore property, it is natural to ask for any relationship between PL and Howe-Moore, and this is an interesting {\it open} question.  In the appendix, Corina Ciobotaru shows that a closed non-compact subgroup of the automorphism group of the $d$-regular tree ($d\geq 3$) that acts 2-transitively on the boundary, satisfies property PL. As a consequence of her result, all {\it known} examples of groups with the Howe-Moore property (see \cite{Cio}) have property PL. 

This paper is a natural continuation of \cite{CTV, CCLTV}, but can be read independently.

\medskip
{\bf Acknowledgements:} Special thanks are due to Yves Cornulier for numerous exchanges and conversations following the joint papers \cite{CTV} and \cite{CCLTV}; in particular he provided Example \ref{Heisenberg} and suggested the use of the map $\phi^x$ in the proof of Theorem \ref{semidir}. We also thank Bachir Bekka for suggesting Proposition \ref{PropBek}, Yves Benoist for sharing his expertise on algebraic groups, and Pierre-Emmanuel Caprace for suggesting Theorem \ref{P.-E.C.} as an improvement of Theorem \ref{algebraicPL}.

\section{Generalities on property $BP_{L^p}$}

The two next results follow immediately from definitions.

\begin{Prop}\label{prop:propPH_triv}
Let $G$ be a locally compact group, and $N$ a closed normal subgroup.
\begin{enumerate}
\item[1)] If $G$ has property $BP_{L^p}$, then so does $G/N$.
\item[2)] If $G$ has property $BP_{L^p}$, and $N$ is not compact, then $G/N$ has property $FL^p$.
\item[3)] If $G/N$ has property $BP_{L^p}$, and $N$ is compact, then $G$ has property $BP_{L^p}$.\hfill$\square$
\end{enumerate}
\end{Prop}

\begin{Prop}\label{cocom} Let $H$ be a closed co-compact subgroup in $G$. If $H$ has property $BP_{L^p}$, then so does $G$. \hfill$\square$
\end{Prop}

\begin{Ex} Let $N=\mathrm{SL}_{2}(\mathbf{R})\times \mathrm{SL}_{2}(\mathbf{R})$, and let $\mathbf{Z}/2\mathbf{Z}$ act on $N$ by exchanging factors. Form the semi-direct product $G=N\rtimes \mathbf{Z}/2\mathbf{Z}$. Clearly $N$ does not have Property $BP_{L^2}$, but $G$ has Property $BP_{L^2}$ by Theorem \ref{algebraic}. This example shows that Property $BP_{L^2}$ is {\it not} inherited by finite index subgroups. 
\end{Ex}

\begin{Ex}\label{SU(n,1)}
    Let $G$ be the universal covering group of
$\mathrm{SU}(n,1)\,(n\geq 1)$. For every $p\geq 1$, the group $G$ does {\it not} have have property
$BP_{L^p}$, by Proposition \ref{prop:propPH_triv} (since the quotient
$G/Z(G)$ of $G$ by the non-compact normal subgroup $Z(G)$, does not
have property $FL^p$). This shows that property $BP_{L^p}$ is not inherited
by non-trivial central extensions.
\end{Ex}

\begin{Rem} There are plenty of discrete groups with
Property $BP_{L^2}$ provided by discrete groups with Property FH. But
we do not know any example of a discrete group with Property $BP_{L^2}$
but without Property FH. It is a result by Peterson-Thom (Theorem 2.6 in \cite{PeTh}) that, if $G$ is a countable group with non-zero first $L^2$-Betti number, containing some infinite amenable subgroup (e.g. $\Z$), then there exists a 1-cocycle with respect to the regular representation, which is neither bounded or proper; so such a group does not have property $BP_{L^2}$, nor the weaker property $BP_0$.
\end{Rem}

Since a locally compact group admitting a proper isometric action on some metric space must be $\sigma$-compact, we have in particular:

\begin{Prop}\label{sigmacon} A group with property $BP_{L^p}$ but without property $FL^p$, is $\sigma$-compact. \hfill$\square$
\end{Prop} 

Recall that a locally compact group $G$ is {\it locally elliptic} if every compact subset is contained in a compact subgroup. Observe that an locally elliptic group is amenable, as a direct limit of compact groups. For an arbitrary locally compact group, the {\it locally elliptic radical} $R_{ell}(G)$ is the unique maximal locally elliptic closed normal subgroup of $G$; see Example 4.D.7(7) in \cite{CorHar}.

\begin{Prop}\label{locfin} Let $G$ be a $\sigma$-compact group with property $BP_{L^p}$. If $G$ is {\it not} compactly generated, then $G$ is locally elliptic and not almost connected. 
\end{Prop}

{\bf Proof:} Let $K$ be a compact set in $G$, and let $U$ be the closed subgroup generated by $K$. Upon replacing $K$ by its union with a compact neighborhood of the identity, we may assume that $U$ is an open subgroup. Let $(K_n)_{n\geq 0}$ be an increasing sequence of compact subsets of $G$, with $K=K_0$ and $G=\bigcup_{n\geq 0}K_n$, and let $U_n$ be the subgroup generated by $K_n$. By assumption $U_n\neq G$. As explained e.g. in the proof of Proposition 2.4.1 of \cite{BHV}, the set $T=\coprod_{n\geq 0}G/U_n$ carries a natural $G$-invariant tree structure such that, for the $G$-representation $\pi$ on $\ell^p$ of the set of oriented edges, there is an unbounded 1-cocycle $b\in Z^1(G,\pi)$. Actually $\|b(g)\|^p=d(gx_0,x_0)$, where $d(.,.)$ is the distance in $T$ and $x_0$ is the trivial coset in $G/U$ (see Proposition 2.3.3 in \cite{BHV}). By property $BP_{L^p}$, this unbounded cocycle $b$ must be proper, in particular vertex stabilizers in $T$ must be compact. So $U=\{g\in G:gx_0=x_0\}$ is compact, i.e. $G$ is locally elliptic. It then follows from Proposition 4.D.3 in \cite{CorHar}, that the connected component of identity in $G$ is compact.
\hfill$\square$

\medskip
Recall from \cite{CMV} that a {\it space with measured walls} is a 4-tuple $(X,\mathcal{W},\mathcal{B},\mu)$ where $\mathcal{W}$ is a set of walls on $X$ (i.e. partitions of $X$ into 2 classes), $\mathcal{B}$ is a $\sigma$-algebra of subsets of $\mathcal{W}$ and $\mu$ is a measure on $\mathcal{B}$ such that, for any $x,y\in X$, the set $\mathcal{W}(x|y)$ of walls separating $x$ from $y$ is in $\mathcal{B}$ and has finite measure.

The kernel $(x,y)\mapsto\mu(\mathcal{W}(x|y))$ is then a pseudo-metric on $X$, called the {\it wall distance}. 

\begin{Prop}\label{measuredwalls} For a locally compact group $G$ and $p\geq 1$, consider the following statements:
\begin{enumerate}
\item[1)] $G$ has property $BP_{L^p}$.
\item[2)] Every action of $G$ on a space with measured walls $X$ is either bounded or proper (when $X$ is endowed with the wall distance).
\end{enumerate}
Then $(1)\Rightarrow(2)$. The converse holds if $1\leq p\leq2$ and $G$ is separable.
\end{Prop}


{\bf Proof:} $(1)\Rightarrow(2)$ follows essentially from the proof of Proposition 3.1 in \cite{CTV} and the remark following it. We recall the main features. If $G$ acts on a space with measured walls $X$, and $x_0$ is some base-point in $X$, there is an affine isometric action $\alpha_X$ of $G$ on $L^p$ of the space of half-spaces in $X$, such that $\|\alpha_X(g)(0)\|_p^p=\mu(\mathcal{W}(gx_0|x_0))$, the measure of the set of walls separating $gx_0$ from $x_0$. So the $G$-action on $X$ is proper (resp. bounded) if and only if $\alpha_X$ is proper (resp. bounded).

$(2)\Rightarrow(1)$. This is a combination of results from \cite{CDH}. Assume $1\leq p\leq2$ and $G$ separable, and let $\alpha$ be an affine isometric action of $G$ on $L^p(\Omega,\nu)$. Fix $v\in L^p(\Omega,\nu)$ and set $\psi(g)=\|\alpha(g)v-v\|_p^p$. In the terminology of Definition 6.5 in \cite{CDH}, the invariant kernel $(g,h)\mapsto\psi(g^{-1}h)$ has type $p$. By Corollary 6.11(1) in \cite{CDH}, the function $\psi$ is conditionally negative definite on $G$ (because $1\leq p\leq 2$). By Theorem 6.25(2) in \cite{CDH}, since $G$ is separable there exists a median space $(X,d)$, a point $x_0\in X$ and a continuous isometric $G$-action on $X$ such that $\sqrt{\psi(g)}=d(gx_0,x_0)$ for every $g\in G$. Finally, by Theorem 5.1 in \cite{CDH}, because $X$ is median it carries a structure of space with measured walls $(X,\mathcal{W},\mathcal{B},\mu)$ such that $d(x,y)=\mu(\mathcal{W}(x|y))$ for every $x,y\in X$, and every isometry of $X$ is an automorphism of $(X,\mathcal{W},\mathcal{B},\mu)$. So $\psi$ is bounded (resp. proper) if and only if the $G$-action on $(X,\mathcal{W},\mathcal{B},\mu)$ is bounded (resp. proper).
\hfill$\square$

\section{Generalities on property PL}

Let $G$ be a locally compact $\sigma$-compact group. Observe that, if $G$ has property PL, then every closed normal subgroup is either compact or co-compact. If $G$ is amenable with property $BP_{L^p}$ for some $p\geq 1$, and $N$ is a closed non-compact normal subgroup, then $G/N$ is both amenable and property $FL^p$, so $G/N$ is compact. So also in this case every closed normal subgroup of $G$ is either compact or co-compact. 

In Theorem E of \cite{CaMo}, Caprace and Monod obtained structural results for compactly generated, locally compact groups $G$ with the property that every non-trivial closed normal subgroup is co-compact. If $G$ is not compact, then $G$ falls into one of the following three cases:
\begin{enumerate}
\item $G$ is isomorphic to a semi-direct product $\R^d\rtimes K$ where $K$ is a compact subgroup of $\mathrm{GL}_d(\R)$ acting irreducibly on $\R^d$;
\item $G$ is a compact extension of a quasi-product of finitely many non-compact, pairwise isomorphic, topologically simple groups, permuted transitively by $Ad(G)$;
\item $G$ is discrete and residually finite.
\end{enumerate}

\begin{Lem}\label{CapMon} Let $G$ be a non-compact, locally compact group. Assume either that $G$ has property PL, or that $G$ is amenable with property $BP_{L^p}$ for some $p\geq 1$.
\begin{enumerate}
\item[a)] If $G$ is not compactly generated, then $G$ is locally elliptic.
\item[b)] If $G$ is compactly generated, then $R_{ell}(G)$ is compact and every non-trivial closed normal subgroup of $G/R_{ell}(G)$ is co-compact.
\end{enumerate}
\end{Lem}

{\bf Proof:} If $G$ is not compactly generated, the result follows from Proposition \ref{locfin} (as PL implies $BP_{L^p}$). If $G$ is compactly generated, then every closed normal subgroup of $G$ is either compact or co-compact: this is obvious if $G$ has property PL; if $G$ is amenable with property $BP_{L^p}$, this follows from the fact that any non-compact compactly generated amenable group admits a proper action on $L^p$. In particular $R_{ell}(G)$ is compact, and $G/R_{ell}(G)$ is a non-compact group, without non-trivial compact normal subgroup, and every non-trivial closed normal subgroup co-compact.  \hfill$\square$

\medskip
For almost connected Lie groups, Lemma \ref{CapMon} cleans things up, as those are compactly generated. An immediate consequence of Lemma \ref{CapMon} and the Caprace-Monod theorem, is:

\begin{Prop}\label{Liealmostcon} Let $G$ be a non-compact almost connected real Lie group. Assume that $G$ either has property PL, or that $G$ is amenable with property $BP_{L^p}$ for some $p\geq 1$. There is a compact normal subgroup $W$ of $G$ such that:
\begin{enumerate}
\item[a)] if $G$ is amenable, then $G/W$ is isomorphic to a semi-direct product $\R^d\rtimes K$ where $K$ is a compact subgroup of $\mathrm{GL}_d(\R)$ acting irreducibly on $\R^d$;
\item[b)] if $G$ is non-amenable, then $G/W$ is a compact extension of a product of finitely many non-compact, pairwise isomorphic, simple Lie groups, permuted transitively by $Ad(G)$.
\end{enumerate}\hfill$\square$
\end{Prop}

\section{Amenable groups}

\subsection{Semi-direct products}
For $A$ a LCA group, we denote by $\hat{A}$ its Pontryagin dual, and by $1_A\in\hat{A}$ the trivial character. Set $\hat{A}^*=\hat{A}\backslash\{1_A\}$.

\begin{Prop}\label{unbounded} Fix $p\geq 1$. Let $A$ be a LCA group, and let $K$ be a compact group of automorphisms. Let $\mu$ be an infinite $K$-invariant Radon measure on $\hat{A}^*$. Assume that, for every compact subset $C\subset A$, we have 
\begin{equation}\label{Lfinite}
\sup_{a\in C}\int_{\hat{A}^*} |\chi(a)-1|^p\,d\mu(\chi)<+\infty.
\end{equation}
For $(a,k)\in A\rtimes K$, set:
$$L(a,k)=\left(\int_{\hat{A}^*} |\chi(a)-1|^p\,d\mu(\chi)\right)^{1/p}.$$
Then $L$ is an unbounded $L^p$-type length function on $A\rtimes K$.
\end{Prop}

{\bf Proof:} Let $\mathcal{F}_\mu$ be the space of $\mu$-measurable functions on $\hat{A}^*$, modulo equality $\mu$-almost everywhere. We define a linear representation $\pi_\mu$ of $A\rtimes K$ on $\mathcal{F}_\mu$ by 
$$(\pi_\mu(a,k)f)(\chi)=\chi(a)f(k^{-1}\cdot\chi)$$
($(a,k)\in A\rtimes K, f\in\mathcal{F}_\mu,\chi\in\hat{A}^*$). Observe that $\pi_\mu$ has no non-zero fixed vector. View the space $L^p(\hat{A}^*,\mu)$ as a subspace of $\mathcal{F}_\mu$: it is invariant under $\pi_\mu$, that induces an isometric representation of $A\rtimes K$ on $L^p(\hat{A}^*,\mu)$. Let $t$ be the translation by the constant function 1 on $\mathcal{F}_\mu$, so that $t(f)=f+1$ for $f\in\mathcal{F}_\mu$. Define an affine action $\alpha_\mu$ of $A\rtimes K$ on $\mathcal{F}_\mu$ by:
$$\alpha_\mu=t^{-1}\circ\pi_\mu\circ t.$$
More precisely, for $(a,k)\in A\rtimes K, f\in\mathcal{F}_\mu,\chi\in\hat{A}^*$:
$$(\alpha_\mu(a,k)f)(\chi)=(\pi_\mu(a,k)f)(\chi)+ \chi(a)-1.$$
Observe that the constant function -1 is the only fixed point of $\alpha_\mu$ in $\mathcal{F}_\mu$.

By assumption $\chi\mapsto \chi(a)-1$ is in $L^p(\hat{A}^*,\mu)$ for every $a\in A$, so $L^p(\hat{A}^*,\mu)$ is $\alpha_\mu$-invariant, and $\alpha_\mu$ defines a continuous affine isometric action of $A\rtimes K$ on $L^p(\hat{A}^*,\mu)$. Then $L(a,k)=\|\alpha_\mu(a,k)(0)\|$ is indeed a $L^p$-type length function on $A\rtimes K$. Since $-1\notin L^p(\hat{A}^*,\mu)$, the action $\alpha_\mu$ has no fixed point in $L^p(\hat{A}^*,\mu)$, so that $L$ is unbounded. Indeed, this follows from the fact that an isometric action on $L^p$ with bounded orbits fixes a point: this is a consequence of the center lemma for $p>1$, and of \cite{BGM} for $p=1$).
\hfill$\square$

\begin{Rem}\label{LfiniteRd} Suppose that $A=\R^d$; every continuous character on $\R^d$ is of the form $x\mapsto\exp(2\pi i\langle x|y\rangle)$, for some $y\in\R^d$ (where $\langle.|.\rangle$ denotes the usual scalar product), so for $(x,k)\in\R^d\rtimes K$ we have
$$L(x,k)^p=\int_{(\R^d)^*} |\exp(2\pi i\langle x|y\rangle)-1|^p\,d\mu(y)= 2^{p/2}\int_{(\R^d)^*}(1-\cos(2\pi\langle x|y\rangle))^{p/2}\,d\mu(y)$$
$$=2^p\int_{(\R^d)^*}|\sin^p(\pi\langle x|y\rangle)|\,d\mu(y).$$
Using $|\sin t|\leq|t|$ and the Cauchy-Schwarz inequality, we see that, to ensure the finiteness condition (\ref{Lfinite}), it is sufficient that $\mu$ has a $p$-th moment, i.e. $\int_{(\R^d)^*}\|y\|^p\,d\mu(y)<+\infty$.
\end{Rem}

We will apply Proposition \ref{unbounded} to semi-direct products $A\rtimes F$, with $F$ finite. In the case $A=\Z$ and $F$ trivial, the next result is due to Edelstein (Theorem 2.1 in \cite{Edel}); for $A=\R$ and $F$ trivial, see Corollary 5.3 in \cite{CTV}.

\begin{Thm}\label{nonBP} Let $A$ be a non-compact, $\sigma$-compact LCA group, and let $F$ be a finite group of automorphisms of $A$. For every $p\geq1$, the semi-direct product $A\rtimes F$ does not have property $BP_{L^p}$.
\end{Thm}

{\bf Proof:} We will use Proposition \ref{unbounded} to construct a specific unbounded $L^p$-type length function on $A\rtimes F$, which will turn out to be not proper. Let $(K_n)_{n>0}$ be a strictly increasing sequence of compact subsets of $A$, with $A=\cup_{n>0} K_n$. Clearly we may assume that each $K_n$ is $F$-invariant. For $a\in A$, let $|a|$ be the unique integer $n>0$ such that $a\in K_n\backslash K_{n-1}$, so that $K_n=\{a\in A: |a|\leq n\}$.

{\bf Claim:} There exists sequences $(a_n)_{n>0}$ in $A$, and $(\chi_n)_{n>0}$ in $\hat{A}^*$ such that:
\begin{itemize}
\item $|a_n|>|a_{n-1}|;$
\item for $f\in F,k\leq n$ we have: $|(f\cdot\chi_k)(a_n)-1|\leq 2^{-n}$;
\item $\max_{|a|\leq |a_{n-1}|} |\chi_n(a)-1|\leq 2^{-n}$ for every $n>0$.
\end{itemize}

Taking the Claim for granted, we define the measure $\mu$ on $\hat{A}^*$ as a sum of Dirac masses:
$$\mu=\sum_{f\in F}\sum_{k>0}\delta_{(f\cdot\chi_k)}.$$
Then $\mu$ is an infinite $F$-invariant Radon measure on $\hat{A}^*$. Moreover for $a\in K_{n}$ we have a uniform bound:
$$\int_{\hat{A}^*} |\chi(a)-1|^p\,d\mu(\chi)=\sum_{f\in F}\sum_{k=1}^{n} |(f\cdot\chi_k)(a)-1|^p + \sum_{f\in F}\sum_{k>n}|(f\cdot\chi_k)(a)-1|^p$$
$$\leq 2^pn|F| + |F|\sum_{k> n}2^{-pk}$$
where the inequality follows from the Claim and $K_n\subset K_{|a_n|}\subset K_{|a_{k-1}|}$ for $k>n$. Then by Proposition \ref{unbounded}:
$$L(a,g)=[\sum_{f\in F}\sum_{k>0}|(f\cdot\chi_k)(a)-1|^p]^{1/p}$$
($(a,g)\in A\rtimes F$) defines an unbounded $L^p$-type length function on $A\rtimes F$. To show that it is not proper, we show that $L$ remains bounded along the unbounded sequence $(a_n,Id_A)_{n>0}$ in $A\rtimes F$. 
But by the Claim:
$$L(a_n,Id_A)^p=\sum_{f\in F}\sum_{k=1}^{n}|(f\cdot\chi_k)(a_n)-1|^p + \sum_{f\in F}\sum_{k>n}|(f\cdot\chi_k)(a_n)-1|^p$$
$$\leq |F|\cdot n\cdot 2^{-pn} + |F|\sum_{k>n}2^{-pk}$$
as $|a_n|\leq|a_{k-1}|$ for $k>n$. So $\lim_{n\rightarrow\infty} L(a_n,Id_A)=0$, so the sequence $(L(a_n,Id_A))_{n>0}$ is bounded.

It remains to prove the Claim. Assume that $a_1,...,a_{n-1}\in A$ and $\chi_1,...,\chi_{n-1}\in\hat{A}^*$ have been constructed. As $A$ is not compact, so that $\hat{A}$ is not discrete, we find $\chi_n\in\hat{A}^*$ such that $\max_{|a|\leq |a_{n-1}|} |\chi_n(a)-1|\leq 2^{-n}$. Let $H$ be the subgroup of $\hat{A}$ generated by the $f\cdot\chi_k$'s, with $f\in F,k\leq n$. Let $\iota$ denote the inclusion of $H$ into $\hat{A}$. We then have the dual homomorphism $\iota^*:\hat{\hat{A}}=A\rightarrow \hat{H}$ with $\hat{H}$ compact. Since $\iota^*$ has dense image (see Corollaire 6 in Chap II.1.7 of \cite{Bour}), and $A$ is not compact, in any complement of a compact set in $A$ we can find $a$ with $\iota^*(a)$ arbitrarily close to the trivial element in $\hat{H}$. In particular we can find $a_n$ with $|a_n|>|a_{n-1}|$ with $|(\iota^*(a_n))(f\cdot\chi_k)-1|\leq 2^{-n}$, i.e. $|(f\cdot\chi_k)(a_n)-1|\leq 2^{-n}$ for $f\in F,k\leq n$.
\hfill$\square$

\begin{Rem} Say that $A=\R^d$ in Theorem \ref{nonBP}, and assume that the finite group $F$ stabilizes some proper, closed, unbounded subgroup $B$ of $\R^d$ (think of $B$ as either a proper linear subspace, or a lattice). Then the proof of Theorem \ref{nonBP} becomes much simpler. Indeed let $y_0\in\R^d$ be a non-zero vector such that the character $x\mapsto\exp(2\pi i\langle x|y_0\rangle)$ is in the annihilator $B^\perp=\{\chi\in\hat{A}:\chi|_B\equiv 1\}$. Set then $\chi_k(x)=\exp(\frac{2\pi i}{k!}\langle x|y_0\rangle)$. Then the measure $\mu=\sum_{f\in F}\sum_{k>0}\delta_{(f\cdot\chi_k)}$ on $(\R^d)^*$ has finite $p$-th moment, so by remark \ref{LfiniteRd} and Proposition \ref{unbounded} the function $L(a,g)=:[\sum_{f\in F}\sum_{k>0}|(f\cdot\chi_k)(a)-1|^p]^{1/p}$ defines an unbounded $L^p$-type length function on $A\rtimes F$. On the other hand, pick any non-zero vector $a_0 \in B$, and set $a_n=n!\cdot a_0$. We claim that $\lim_{n\rightarrow\infty} L(a_n,Id_A)=0$. Indeed observe that $(f\cdot\chi_k)(a_n)=1$ for $f\in F, k\leq n$, so that as in remark \ref{LfiniteRd}:
$$L(a_n,Id_A)^p=\sum_{f\in F}\sum_{k>n} |(f\cdot\chi_k)(a_n)-1|^p=2^p\sum_{f\in F}\sum_{k>n}|\sin^p(\pi\frac{n!}{k!}\langle a_0|y_0\rangle)|$$
$$\leq (2\pi)^p\|a_0\|^p\|y_0\|^p |F|\sum_{k>n} (\frac{n!}{k!})^p\leq  (2\pi)^p\|a_0\|^p\|y_0\|^p |F|\sum_{k>n} (\frac{1}{n})^{p(n-k)}$$
using the bound $\frac{n!}{k!}<\frac{1}{n^{n-k}}$ for $k>n$. So we have 
$$L(a_n,Id_A)^p\leq (2\pi)^p\|a_0\|^p\|y_0\|^p |F|\cdot\frac{1}{n^p -1},$$ 
establishing the claim. The choice of the weights in defining $\mu$ and the sequence $(a_n)_{n>0}$, is inspired by the proof of Theorem 2.1 in \cite{Edel}. We will come back to finite groups stabilizing a lattice in $\R^d$, in Corollary \ref{motion} below.
\end{Rem}

\medskip 
Here is a noteworthy consequence of Theorem \ref{nonBP}.

\begin{Cor}\label{crystal}
Let $\Gamma$ be an infinite, finitely generated, virtually abelian group. Then  $\Gamma$ does not have property $BP_{L^p}$, for every $p\geq1$.
\end{Cor}

{\bf Proof:}  
Write $\Gamma$ as the central term of a short exact sequence 
\begin{equation}\label{short}
0\rightarrow\Z^n\rightarrow\Gamma\stackrel{p}{\rightarrow} F\rightarrow 1
\end{equation}
with $F$ finite. We claim that $\Gamma$ embeds as a co-compact lattice in a semi-direct product $G=\R^n\rtimes F$. Since $G$ does not have property $BP_{L^p}$ (by Theorem \ref{nonBP}), the Corollary follows from Proposition \ref{cocom}.

To prove the claim, let $c\in Z^2(F,\Z^n)$ be the 2-cocycle on $F$ describing the extension (\ref{short}). The group $\Z^n$ becomes an $F$-module through the conjugation action of $F$, and $c(g,g)=s(g)s(g')s(gg')^{-1}$ for some section $s:F\rightarrow\Gamma$ of the map $p$. Now the $F$-action on $\Z^n$ canonically extends to $\R^n$, and the law 
$$(v,g)(v',g')=(v+g\cdot v'+c(g,g'),gg')\;(v\in\R^n,g\in F)$$
defines on $G=\R^n\times F$ the structure of an almost connected Lie group in which $\Gamma$ embeds as a co-compact lattice. 
Since $H^2(F;\R^n)=0$, the extension
$$0\rightarrow\R^n\rightarrow G \rightarrow F\rightarrow 1$$
splits, so that $G=\R^n\rtimes F$. This proof was inspired by the proof of Theorem 1 in \cite{AuKu}.
\hfill$\square$

\subsection{Abelian groups}

The following Lemma can be deduced from the proof of Proposition 2.5.9 in \cite{BHV} (where it is proved for solvable groups and $p=2$). We include the simple proof for locally compact abelian (LCA) groups.

\begin{Lem}\label{abelianFLp} Fix $p\geq 1$. A LCA group $A$ has property $FL^p$ if and only if $A$ is compact.
\end{Lem}

{\bf Proof:} One implication is trivial. For the non-trivial one, let $A$ be a LCA group with property $FL^p$. We consider two cases:
\begin{itemize}
\item $A$ is discrete. Assume by contradiction that $A$ is infinite. Since every infinite abelian group has a countably infinite quotient (see Theorem 2.5.2 in \cite{Rudin}), we may assume that $A$ is countably infinite with property $FL^p$. But a countable group with property $FL^p$ is finitely generated (by Corollary 2.4.2 in \cite{BHV}), hence $A$ is isomorphic to $\Z^n\oplus F$, with $n>0$ and $F$ finite abelian. But such a group does not have property $FL^p$, so a contradiction is reached.
\item $A$ is arbitrary. By structure theory for LCA groups (see Theorem 2.4.1 in \cite{Rudin}), $A$ admits an open subgroup $U$ of the form $U=K\times\R^m$, for some $m\geq 0$ and some compact group $K$. The group $A/U$ is discrete with property $FL^p$, so it is finite by the first case of the proof, i.e. $U$ has finite index in $A$. So $U$ has property $FL^p$ too (by Proposition 2.5.7 in \cite{BHV}). This clearly forces $m=0$, so $A$ is compact.
\hfill$\square$
\end{itemize}

\bigskip
{\bf Proof of Theorem \ref{abelian}:} Implications $(c)\Rightarrow(a)\Rightarrow(b)$ are clear. We prove $(b)\Rightarrow(c)$ by contradiction. So suppose there is a non-compact LCA group $A$ with property $BP_{L^p}$, for some $p\geq1$. As $A$ is not compact, $A$ does not have property $FL^p$, by Lemma \ref{abelianFLp}. Since $A$ has property $BP_{L^p}$, the group $A$ must be $\sigma$-compact, by Proposition \ref{sigmacon}. But this contradicts Theorem \ref{nonBP}.
\hfill$\square$ 

\medskip
From Theorem \ref{abelian} we deduce immediately:

\begin{Cor}\label{abelianized} Let $G$ be a locally compact group with property $BP_{L^p}$, for some $p\geq 1$. Then $G/\overline{[G,G]}$ is compact.
\hfill$\square$
\end{Cor}

\subsection{Centers: proof of Theorem \ref{centralext}}
Let $S$ be a compact generating subset of $G$.
Since $G$ does not have property $FL^p$, it admits an affine isometric action $\sigma$ on some $L^p$-space $E$ with non-zero displacement: 
$$\inf_{x\in E} \sup_{s\in S}\|\sigma(s)x-x\|>0.$$
Indeed, this follows from \cite[3.8.D]{Gro}, taking the ``energy" $E(x)$ to be $\sup_{s\in S}\|\sigma(s)x-x\|$. In other words, letting $\pi$ and $b$ be respectively the linear part and cocycle part of $\sigma$, we have that $b$ is non-trivial in reduced first cohomology. Let $Z(G)$ denote as usual the center of $G$.
By \cite[Proposition 2.6]{BFGM}, we have a $\pi(G)$-invariant continuous decomposition $E=E_1\oplus E_2\oplus E_3$, where $E_1$ is the space of $\pi(G)$-invariant vectors, and $E_1\oplus E_2$ the space of $Z(G)$-invariant vectors. The projection $b_1$ of $b$ on $E_1$ is a group homomorphism $G\rightarrow E_1$, which by Corollary \ref{abelianized} is zero. Observe then that the projection $b_2$ of $b$ on $E_2$ vanishes on $Z(G)$. Indeed, the cocycle relation shows that $b_2(z)$ is a $\pi(G)$-invariant vector for all $z\in Z(G)$. So assuming by contradiction that $Z(G)$ is not compact, by property $BP_{L^p}$ the affine action corresponding to $b_2$ is bounded, so it has a fixed point and $b_2$ is a coboundary.
Finally, as the center $Z(G)$ is non-compact, the projection $b_3$ is trivial in $\overline{H}^1(\pi,E)$ by \cite[Corollary 5]{BRS}, finally implying that $b$ itself is an almost coboundary: this is a contradiction.
\hfill$\square$

\medskip
Note that the above proof really needs $p>1$ to appeal to Proposition 2.6 of \cite{BFGM} and to Corollary 5 of \cite{BRS}.

\subsection{Proof of Theorem \ref{semidir}}

\begin{Lem}\label{Liealgebra} Let $F$ be a local field of characteristic 0, $V=F^d$, and $K$ an infinite compact subgroup of $\mathrm{GL}(V)$, acting irreducibly on $V$. Denote by $\mathfrak{k}$ the Lie algebra of $K$. For every non-zero $x\in V$, there exists $X\in\mathfrak{k}$ such that $Xx\neq 0$.
\end{Lem}

{\bf Proof:} Contraposing, we assume that there is a non-zero vector $x\in V$ such that $Xx=0$ for every $X\in\mathfrak{k}$, and will show that $K$ is finite. Let $W$ be the space of vectors $v\in V$ such that $Xv=0$ for every $X\in\mathfrak{k}$: this is clearly a $K$-invariant subspace, so by irreducibility we have $W=V$. This implies $\mathfrak{k}=0$ and hence $K$ is finite.
\hfill$\square$

\medskip
The next lemma says that, if $G$ is as in Theorem \ref{semidir}, it is close to being a semi-direct product.

\begin{Lem}\label{almostsplit} Let $G$ be as in Theorem \ref{semidir}. There exists a compact subgroup $C$ of $G$ such that $G=VC$.
\end{Lem}

{\bf Proof:} If $F=\R$, it is a classical fact (see Theorem 2.3 in Chapter III of \cite{Hoch}) that any extension of a finite-dimensional real vector space by a compact group, is split. For $F$ non-Archimedean, we observe that $V$ is locally elliptic and appeal to the fact that local ellipticity is preserved by extensions (see Proposition 4.D.4(2) in \cite{CorHar}): hence $G$ is locally elliptic. So if $K$ is a compact set that surjects onto $G/V$ by the quotient map $G\rightarrow G/V$, the set $K$ is contained in a compact subgroup $C$ of $G$, and clearly $G=VC$.
 \hfill$\square$

\begin{Ex}\label{Heisenberg} Let $G$ be the following closed subgroup of the Heisenberg group over $\Q_p$:
$$G=\left\{\left(\begin{array}{ccc}1 & x & z \\0 & 1 & y \\0 & 0 & 1\end{array}\right):x,y\in\Z_p,z\in\Q_p\right\}.$$
Then $G$ is a central extension of $V=\Q_p$ by $\Z_p^2$. The extension is not split as $G$ is not abelian. However we have $G=VC$ where $C$ is the Heisenberg group over $\Z_p$.
\end{Ex}

{\bf Proof of Theorem \ref{semidir}:} The implication $(a)\Rightarrow (b)$ is obvious.

$(b)\Rightarrow(c)$ Assume that $G$ has property $BP_{L^p}$. By Proposition \ref{measuredwalls} we may assume $p>1$. By Lemma \ref{almostsplit}, we can write $G=VC$ for some compact subgroup $C$ of $G$. Observe that $V\cap C$ is normal in $G$. Denote by $\alpha:V\rightarrow V/(V\cap C)$ and $\beta: C\rightarrow C/(V\cap C)$ the quotient maps. Then the map 
$$G\rightarrow V/(V\cap C)\rtimes C/(V\cap C):g=vc\mapsto (\alpha(v),\beta(c))$$
is well-defined and is a continuous surjective homomorphism. So the semi-direct product $V/(V\cap C)\rtimes C/(V\cap C)$ has property $BP_{L^p}$ and Theorem \ref{nonBP} implies that $G/V\simeq C/(V\cap C)$ is infinite. If $W$ is a non-zero $G/V$-invariant linear subspace in $V$, then $W$ is a normal subgroup in $G$. By Proposition \ref{prop:propPH_triv}, the quotient $G/W$ has property $FL^p$, so it is compact as $G/W$ is also amenable. Hence $W=V$, i.e. $G/V$ acts irreducibly on $V$.

$(c)\Rightarrow (a)$ Set $K=G/V$. Assume that $K$ acts irreducibly on $V$ and is infinite, so that its Lie algebra $\mathfrak{k}$ is non-zero. We proceed in several steps. 
\begin{itemize}
\item For $x\in V\backslash\{0\}$, set $\phi^x:K^d\rightarrow V:(g_1,...,g_d)\mapsto g_1x+...+g_dx$. We claim that the image of $\phi^x$ contains some open set in $V$. For this, it is enough to show that the differential $D\phi^x_{(g_1,...,g_d)}$ has rank $d$ for some $(g_1,...,g_d)\in K^d$. But for $(X_1,...,X_d)\in\mathfrak{k}^d$ we have:
$$D\phi^x_{(g_1,...,g_d)}(X_1,...,X_d)=g_1X_1x+...+g_dX_dx.$$
By Lemma \ref{Liealgebra}, we find $X\in\mathfrak{k}$ such that $Xx\neq 0$. As $Xx$ is a cyclic vector for $K$ (because $K$ acts irreducibly), we find $g_1,...,g_d\in K$ such that $\{g_1Xx,...,g_dXx\}$ is a basis of $V$. This means that $D\phi^x_{(g_1,...,g_d)}$ has rank $d$.
\item Endow $V=F^d$ with the $\ell^\infty$ norm $\|x\|=\max_{1\leq i\leq d}|x_i|$. For $x\in V\backslash\{0\}$, set $\psi^x:K^{2d}\rightarrow V:(g_1,...,g_d,h_1,...,h_d)\mapsto \phi^x(g_1,...,g_d)+\phi^{-x}(h_1,...,h_d)$. By the previous point, the image of $\psi^x$ contains some open set around 0. Let $\epsilon_x$ denote the radius of the largest open ball centered at 0 and contained in the image of $\psi^x$. Since $\psi^x$ depends smoothly on $x$, the function $x\mapsto\epsilon_x$ on the unit sphere of $V$, is bounded below by some positive $\epsilon>0$. 
\item As in Lemma \ref{almostsplit}, write $G=VC$. Let $L$ be a length function on $G$, let $N>0$ be such that $L|_C\leq N$. For $x\in V, g\in C$, the relation $gxg^{-1}=g(x)$ in $G$ implies: $$L(g(x))\leq L(x) + 2N.$$
 Assume that $L$ is not proper, so that there is a sequence $(x_n)_{n>0}$ in $V$, with $\|x_n\|\rightarrow\infty$, and a constant $M>0$ such that $L(x_n)\leq M$ for every $n>0$. Fix $y\in V$, and choose $n$ large enough so that $\|y\|<\epsilon\|x_n\|$. This implies that $y$ is in the image of $\psi^{x_n}$, say $y=\sum_{i=1}^d g_i(x _n)+\sum_{i=1}^dh_i(-x_n)$ for suitable $g_1,...,g_d,h_1,...,h_d\in K$. Then
 $$L(y)\leq \sum_{i=1}^d L(g_i(x_n))+\sum_{i=1}^d L(h_i(-x_n))\leq 2d(L(x_n)+2N)$$
 $$\leq 2d(M+2N).$$
 Finally for $y\in V, c\in C$ we have: $L(yc)\leq L(y)+L(c)\leq 2d(M+2N)+N$, meaning that $L$ is bounded.
\end{itemize}\hfill$\square$

\medskip

\subsection{Proof of Theorem \ref{algebraicamenable}}

The implication $(a)\Rightarrow (b)$ is trivial while $(c)\Rightarrow (a)$ follows from Theorem \ref{semidir} together with the fact that property PL is stable under extensions with compact kernels (see Lemma 3.1 in \cite{Cor}).

To prove $(b)\Rightarrow(c)$, let $G$ be non-compact amenable in $\mathcal{G}_{LA}$, with property $BP_{L^p}$. We will repeatedly appeal to what Proposition \ref{prop:propPH_triv} says for amenable groups: if a locally compact amenable group has property $BP_{L^p}$, then every closed normal subgroup is either compact or co-compact. The case of almost connected Lie groups follows from Proposition \ref{Liealmostcon}(a) and Theorem \ref{nonBP}, so we may focus on the algebraic case, i.e. $G=\mathbb{G}(F)$ with $F$ a non-Archimedean local field $F$ of characteristic 0. Let $\mathcal{O}$ be the valuation ring of $F$ and $\pi$ be a uniformizer, so that $F^\times=\mathcal{O}^\times \pi^\Z.$

Let $G^0$ be the Zariski-connected component of identity, and $R_u(G^0)$ its unipotent radical. 
We proceed in two steps:
\begin{itemize}
\item We claim that the unipotent radical $R_u(G^0)$ is non-trivial. Suppose by contradiction that it is trivial, i.e. $G^0$ is reductive. Consider the Levi decomposition $G^0=R(G^0)S$; as $G^0$ is non-compact amenable, $S$ is compact/anisotropic and the radical $R(G^0)$ is a non-compact torus, say $R(G^0)\simeq (F^\times)^r$ with $r>0$. 
Let $T\simeq(\mathcal{O}^\times)^r$ be the unique maximal compact subgroup of $R(G^0)$: then $TS$ is the unique maximal compact subgroup of $G^0$, so it is normal in $G$. The quotient $G/TS$ contains $\Z^r$ with finite index. Because of the assumption $G/TS$ has property $BP_{L^p}$, contradicting Corollary \ref{crystal}. 
\item $[R_u(G^0),R_u(G^0)]=\{1\}$, i.e. $R_u(G^0)$ is abelian (otherwise $[R_u(G^0),R^(G^0)]$ is a non-compact and non-co-compact closed normal subgroup in $G$). So $R_u(G^0)\simeq F^d$ for some $d\geq 1$. Since $R_u(G^0)$ is normal in $G$, it is co-compact. The result then follows from $(b)\Rightarrow(c)$ in Theorem \ref{semidir}.
\end{itemize}

\section{Non-amenable groups: proof of Theorem \ref{algebraicPL}}

It is actually possible to weaken the assumption of Theorem \ref{algebraicPL} rather drastically. For this we need two more definitions.

A locally compact group is {\it locally linear} if it admits an open subgroup which is linear over some local field. We also define the class of {\it elementary groups} as the smallest class of locally compact totally disconnected groups containing all discrete groups, all profinite groups, and closed under group extensions and directed unions of open subgroups. 

\begin{Thm}\label{P.-E.C.} Let $G$ be a non-amenable locally compact group. Assume moreover that $G$ is either an almost connected Lie group, or a non-elementary locally linear totally disconnected group. The following are equivalent:
\begin{enumerate}
\item[a)] $G$ has property PL;
\item[b)] every closed normal subgroup of $G$ is either compact or co-compact;
\item[c)] there exists a compact normal subgroup $W$ and a closed co-compact normal subgroup $N$ of $G$, with $W\triangleleft N$, such that $N/W$ is isomorphic to a direct product $S_1\times...\times S_n$ of simple algebraic groups over some local field $F$, and the simple factors $S_1,...,S_n$ are permuted transitively under $Ad(G)$.
\end{enumerate}
\end{Thm}

In particular Theorem \ref{P.-E.C.} applies to any non-amenable group of the form $\mathbb{G}(F)$, the group of $F$-rational points of a linear algebraic group $\mathbb{G}$ defined over a non-Archimedean local field $F$ of {\it any} characteristic.

\medskip
{\bf Proof of Theorem \ref{P.-E.C.}:}

$(c)\Rightarrow(a)$ Let $G, W$ be as in c). By Lemma 3.1 in \cite{Cor}, it is enough to show that $G/W$ has property PL. So let $L$ be an unbounded length function on $G/W$, we show that $L$ is proper. As $M:=N/W$ is co-compact, $L|_M$ is unbounded. So there exists some index $i$ such that $L|_{S_i}$ is unbounded. Say $i=1$. By assumption, for $j=2,...,r$, there exists $g_j\in G$ such that $Ad(g_j)(S_j)=S_1$. Then, for $s_j\in S_j$ we have by the triangle inequality: $L(g_js_jg_j^{-1})\leq L(s_j)+2L(g_j)$, so that $L|_{S_j}$ is unbounded too. By Theorem 1.4 in \cite{Cor}, $L|_{S_i}$ is proper for every $i=1,...,r$. By Lemma 1.7 in \cite{Cor}, $L|_M$ is proper. So $L$ is proper.

$(a)\Rightarrow(b)$ We already observed that, in a locally compact $\sigma$-compact group with property PL, every closed normal subgroup is either compact or co-compact.

$(b)\Rightarrow(c)$ The Lie group case follows immediately from the already quoted Caprace-Monod theorem (Theorem E in \cite{CaMo}) and the discussion preceding Proposition \ref{Liealmostcon}. 

For the totally disconnected case, we appeal to a structural result by Caprace and Stulemeijer (Corollary 1.2 in \cite{CaStu}): if $G$ is totally disconnected and locally linear, there exists closed characteristic subgroups $W\triangleleft N\triangleleft G$ such that $W$ is elementary, $N/W$ (if non-trivial) is a product $S_1\times...\times S_n$ of topologically simple algebraic groups over local fields $F_1,..., F_n$ (in particular $S_i$ is compactly generated and abstractly simple), and $G/N$ is elementary. In view of the assumption that $G$ is non-elementary, $N/W$ is non-trivial, hence non-compact, in our case. If we assume that all closed normal subgroups of $G$ are either compact or co-compact, we get that $W$ is compact and $N$ is co-compact. Finally $Ad(G)$ acts transitively on the simple factors of $N/W$, since a proper orbit would allow to construct a closed normal subgroup of $G$ that is neither compact nor co-compact.
\hfill$\square$



\section{Property $BP_{L^2}$ in particular}

\subsection{Proof of Theorem \ref{algebraic}}

The implication $(b)\Rightarrow(a)$ is clear: if $G$ has the described form, then by Theorem \ref{algebraicPL} the group $G$ has property PL, {\it a fortiori} it has property $BP_{L^2}$.

The proof of $(a)\Rightarrow(b)$ is very similar in spirit to the proof of $(b)\Rightarrow (c)$ in Theorem \ref{algebraicPL}. Let $G$ be either an almost connected Lie group, or $G=\mathbb{G}(F)$, a linear algebraic group over a local field $F$ of characteristic 0. Let $\mathfrak{g}$ be the Lie algebra of $G$, let $\mathfrak{g}=\mathfrak{r}\rtimes\mathfrak{s}$ be a Levi decomposition, write $\mathfrak{s}=\mathfrak{s}_c\oplus \mathfrak{s}_{nc}$, where $\mathfrak{s}_c$ stands for the compact/anisotropic factors, and $\mathfrak{s}_{nc}$ stands for the non-compact/isotropic factors. Assume $G$ non-amenable, so that $\mathfrak{s}_{nc}\neq 0$.

Suppose that $G$ has property $BP_{L^2}$ but not FH. Then $G$ admits a proper isometric action on a Hilbert space, i.e. $G$ has the Haagerup property. By Theorem 1.10 in \cite{CoJLT}, this implies $[\mathfrak{r},\mathfrak{s}_{nc}]=0$, and all simple factors of $\mathfrak{s}_{nc}$ have $F$-rank 1, and are locally isomorphic to $\mathfrak{so}(n,1)$ or $\mathfrak{su}(n,1)\;(n\geq 2)$ if $F=\R$ or if $G$ is Lie and almost connected. By property $BP_{L^2}$, any quotient of $G$ by a closed non-compact normal subgroup, must have property FH (see Proposition \ref{prop:propPH_triv}). We now distinguish the two cases.

\subsubsection{The algebraic case}

Let $G^0$ be the Zariski-connected component of identity of $G$.   
The radical $R(G^0)$ is compact. Suppose not: as $R(G^0)$ is characteristic in $G^0$, it is a closed non-compact subgroup of $G$, and the quotient $G/R(G^0)$ does not have property FH, contradicting property $BP_{L^2}$ of $G$. 

Let $S_c$ (resp. $S_{nc}$) be the Zariski-connected subgroup of $G^0$ corresponding to $\mathfrak{s}_c$ (resp. $\mathfrak{s}_{nc}$). Since $\mathfrak{s}_{nc}$ is a characteristic ideal in $\mathfrak{g}$, the subgroup $S_{nc}$ is characteristic in $G^0$, hence normal and co-compact in $G$. Finally $Ad(G)$ acts transitively on the simple factors of $S_{nc}$, since a proper orbit would allow to construct a quotient of $G$ by a closed non-compact normal subgroup, not having property FH. Let $W=Z(S_{nc})$ be the center of $S_{nc}$; then the subgroup $N:=S_{nc}/W$ is the desired subgroup of $G/W$.

\subsubsection{The Lie case}

Let $G$ be a non-amenable, almost connected Lie group with $BP_{L^2}$ and without FH. Let $G^0=RS$ be a Levi decomposition of the connected component of identity $G^0$ of $G$. Then $R$ is compact (otherwise, as above, there is a quotient by a non-compact normal subgroup, not having FH), so $R$ is a compact torus, and $G^0$ is reductive. 

Let $S_c$ and $S_{nc}$ be the analytic subgroups corresponding to $\mathfrak{s}_c$ and $\mathfrak{s}_{nc}$ respectively. Note that $S_{nc}$ is closed in $G$ because $R$ is compact. The subgroup $S_{nc}$ is characteristic in $G^0$ hence normal in $G$, so $S_{nc}$ is co-compact in $G$. As above, $Ad(G)$ permutes the simple factors of $S_{nc}$ transitively.

Set $W:=Z(S_{nc})$, the center of $S_{nc}$. Since $W$ is normal in $G$ and $G/W$ does not have property FH, the subgroup $W$ must be finite. As in the algebraic case, we set $N:=S_{nc}/W$, and it is the desired subgroup of $G/W$.

\subsection{Link with the Howe-Moore property}

Let $H$ be a closed subgroup of the locally compact group $G$. We recall Definition 1.3 of \cite{CCLTV}:

\begin{Def} The pair $(G,H)$ has the relative Howe-Moore property if every unitary representation $\pi$ of $G$, either has $H$-invariant vectors, or is such that $\pi|_H$ is a $C_0$-representation. The group $G$ is a Howe-Moore group if the pair $(G,G)$ has the relative Howe-Moore property. 
\end{Def}

\begin{Prop}\label{PropBek} Let $N$ be a closed, co-compact normal subgroup of $G$. If the pair $(G,N)$ has the relative Howe-Moore property, then $G$ has property $BP_{L^2}$. In particular every Howe-Moore group has property $BP_{L^2}$.
\end{Prop}

{\bf Proof:} Let $\pi$ be a unitary representation of $G$, and let $b$ be a 1-cocycle with respect to $\pi$. Set $\psi(g)=\|b(g)\|^2$.
Assuming that $(G,N)$ has the Howe-Moore property, we must prove that $\psi|_N$ is either bounded or proper. 
So suppose that $\psi|_N$ is unbounded.

By Sch\"onberg's theorem, for $t>0$, the function $\phi_t(g)=e^{-t\psi(g)}$ is positive definite on $G$. So there exists a Hilbert space
$\HH_t$ and a unitary representation $\pi_t$ of $G$ on $\HH_t$, with a cyclic vector $\xi_t\in\HH_t$, such that:
$$\phi_t(g)\,=\,\langle\pi_t(g)\xi_t|\xi_t\rangle$$
for every $g\in G$.
\medskip

{\bf Claim:} $\pi_t$ has no non-zero $N$-fixed vector.

To see this, let $\HH^0$ be the space of $N$-fixed vectors in $\HH_t$, and $\HH^{\perp}$ be its orthogonal complement.  We must show that
$\HH^0=0$. Observe that $\HH^0$ and $\HH^{\perp}$ are $G$-invariant, as $N$ is normal in $G$. For $\xi\in\HH$, write $\xi=\xi^0 +\xi^{\perp}$
in the decomposition $\HH_t =\HH^0\oplus\HH^{\perp}$. As $\xi_t$ is cyclic, it is enough to show that $\xi^0_t =0$. But, for $h\in N$:
$$\phi_t(h)\,=\,\langle\pi_t(h)\xi_t^{\perp}|\xi_t^{\perp}\rangle + \|\xi_t^0\|^2.$$
As $\psi|_N$ is unbounded, we can find a sequence $(h_n)_{n>0}$ in
$N$ such that $\lim_{n\rightarrow\infty}\psi(h_n)=+\infty$. Then
$\lim_{n\rightarrow\infty}\phi_t(h_n)=0$. On the other hand,
coefficients of $\pi_t|_N$ on $\HH^{\perp}$ are $C_0$, by the
Howe-Moore property for $(G,N)$. So $\lim_{n\rightarrow\infty}
\langle\pi_t(h_n)\xi_t^{\perp}|\xi_t^{\perp}\rangle\,=\,0$. Hence
$\|\xi_t^0 \|=0$, proving the claim.
\medskip

From the claim, plus the fact that $(G,N)$ has the Howe-Moore
property, we deduce that $\phi_t|_N$ is a $C_0$-function. This is
equivalent to saying that $\psi|_N$ is proper. \hfill$\square$

\medskip
We will see in Example \ref{euclidean} below that the converse of Proposition \ref{PropBek} does not hold in general.

We revisit semi-direct products of the form $V\rtimes K$, where $V=\R^d\,(d\geq 2)$ and $K$ is a closed subgroup of the unitary group $U(d)$. Denote by $K^0$ the connected component of identity of $K$.

\begin{Cor}\label{motion} Consider the following statements:
\begin{enumerate}
\item[a)]  $K^0$ acts irreducibly on $V$.
\item[a')] The pair $(V\rtimes K,V)$ has the relative Howe-Moore property.
\item[b)] $K$ is infinite and acts irreducibly on $V$.
\item[b')] $V\rtimes K$ has property $BP_{L^2}$.
\item[c)] $K$ stabilizes no proper closed unbounded subgroup of $V$.
\item[d)] $K$ acts irreducibly on $V$.
\end{enumerate}
Then $(a)\Leftrightarrow(a')\Rightarrow(b)\Leftrightarrow(b')\Rightarrow(c)\Rightarrow(d)$. If $K$ is connected, all those statements are equivalent.
\end{Cor}

{\bf Proof:} $(a)\Leftrightarrow(a')$ follows immediately from Theorem 4.5 in \cite{CCLTV}.

$(a)\Rightarrow(b)$ follows by observing that, if $K^0$ acts irreducibly, then $K^0$ is non-trivial, hence infinite.

$(b)\Leftrightarrow(b')$ is Theorem \ref{semidir} above.

$(b')\Rightarrow(c)$ If the semi-direct product $V\rtimes K$ has property $BP_{L^2}$ then every closed normal subgroup of $V\rtimes K$ is either compact or co-compact. This rules out any proper closed unbounded $K$-invariant subgroup of $W$.

$(c)\Rightarrow (d)$ is trivial, and so is $(d)\Rightarrow(a)$ when $K=K^0$. \hfill$\square$

\medskip
In general the implications $(a)\Rightarrow(b),\,(b)\Rightarrow(c),\,(c)\Rightarrow(d)$ cannot be reversed, as the following examples show. 

\begin{Ex}\label{euclidean}\begin{enumerate}
\item Let $K$ be the semi-direct product $K=(\mathrm{SO}(2)\times \mathrm{SO}(2))\rtimes C_2$, where $C_2$ acts by flipping the two factors. Then $K$ acts on $\R^4=\R^2\oplus\R^2$, with the first (resp. second) copy of $\mathrm{SO}(2)$ acting by rotations on the first (resp. second) copy of $\R^2$, and $C_2$ flipping the two copies of $\R^2$. Then $K$ acts irreducibly on $\R^4$ but $K^0=\mathrm{SO}(2)\times \mathrm{SO}(2)$ acts reducibly. So $(b)\Rightarrow(a)$ does not hold in general. Since $\R^4\rtimes K$ has property $BP_{L^2}$ but the pair $(\R^4\rtimes K,\R^4)$ does not have the relative Howe-Moore property, the same example shows that the converse of Proposition \ref{PropBek} fails in general. 
\item Let $C_n$ denote the cyclic group of order $n\geq 2$. Let $C_n$ act on $\R^2$ by rotations of angles a multiple of $2\pi/n$. For $n\neq 2,3,4,6$, the group $K=C_n$ stabilizes no proper closed unbounded subgroup of $\R^2$. So $(c)\Rightarrow(b)$ does not hold in general.
\item Consider the same action of $C_n$ by rotations on $\R^2$, this time with $n=3,4,6$: the action is irreducible but stabilizes a lattice in $\R^2$. So $(d)\Rightarrow(c)$ does not hold in general.
\end{enumerate}
\end{Ex}

\section{Groups acting on trees with property PL\\ \small{Appendix by Corina CIOBOTARU}}

In this appendix we prove that closed non-compact subgroups of $\mathrm{Aut}(T_d)^{+}$ that act $2$-transitively on $\partial T_d$ have property $PL$. Beside linear examples as $\mathrm{SL}(2, \Q_p)$, the latter family of groups contains examples of non-compact locally compact groups that are non-linear, at least in characteristic 0: those are the universal groups $U(F)^{+}$ introduced by Burger--Mozes in \cite[Section~3]{BM00a}.

\medskip
We denote by $T_d$ a $d$-regular tree, with $d \geq 3$, and by $\mathrm{Aut}(T_d)$ its group of automorphisms, which is a locally compact group with respect to the compact-open topology. Let $\mathrm{Aut}(T_d)^{+}$ be the group of all type-preserving automorphisms of $T_d$. By a \textbf{type-preserving automorphism} of $T_d$ we mean one that preserves an orientation of $T_d$ that is fixed in advance; this is the same as saying that the automorphism acts without inversion. We denote by $\partial T_d$ the set of endpoints of $T_d$ (they are also called the ideal points of $T_d$) and we call $\partial T_d$ the boundary of $T_d$.  For every two points $x,y \in T_d \cup \partial T_d$ we denote by $[x,y]$ the unique geodesic between $x$ and $y$ in $T_d \cup \partial T_d$. 

For $G \leq \mathrm{Aut}(T_d)$ and $x, y \in T_d \cup \partial T_d$ we define $$G_{[x,y]}:=\{ g \in G \; \vert \; g \text{ fixes pointwise the geodesic } [x,y]\}.$$ In particular, $G_{x}= \{g \in G \; \vert \; g(x)=x\}$. For $\xi \in \partial T_d$ we define $G_{\xi}:= \{g \in G \; \vert \; g(\xi)=\xi \}$ and $G_{\xi}^{0}:= \{g \in G \; \vert \; g(\xi)=\xi \text{ and } g \text{ fixes at least one vertex of } T_d\}$. Notice that $G_{\xi}$ can contain hyperbolic elements; if this is the case then $G_{\xi}^{0} \lneq G_{\xi}$.

\medskip
For the remaining of the appendix we consider $G$ to be a closed non-compact subgroup of $\mathrm{Aut}(T_d)^{+}$ that acts $2$-transitively on $\partial T_d$. One easily sees \cite{Ti70}  that $G$ contains at least one hyperbolic element. Typical examples of such subgroups $G$ are $\mathrm{SL}(2, \Q_p)$ and the universal groups introduced by Burger--Mozes in \cite[Section 3]{BM00a}. Those groups are moreover topologically simple. We will see below that the universal groups are not linear.

\begin{Def}
\label{def::U_alpha}
Let $a$ be a hyperbolic element of $G$. Corresponding to $a$ we define the set
\[ 
U^{+}_{a}:= \{ g \in G \; | \; \lim\limits_{n \to \infty} a^{-n}g a^{n}=e \}. 
\]

Notice that $U^{+}_{a}$ is a subgroup of $G$. It is called the \textbf{contraction group} corresponding to $a$, and in general it is not a closed subgroup of $G$.  In the same way, but using $a^{n}ga^{-n}$ we define $U^{-}_{a}$. For example $U_{a}^{\pm}$ are closed when $G=\mathrm{SL}(2, \Q_p)$ and not closed when $G$ is the universal group $U(F)^{+}$ of Burger--Mozes.
\end{Def}

Let us recall some important properties of $G$, when $G$ is $2$--transitive on $\partial T_d$. Let us fix for what follows a hyperbolic element $a$ of $G$ and denote by $(\xi_{-}, \xi_{+}) \subset \T$ the translation axis of $a$, where $\xi_{-}, \xi_{+} \in \partial T_d$ are the repelling and respectively, the attracting endpoints of $a$. Without loss of generality, we can assume from now on that $a$ is of translation length $2$ (see \cite[Example 4.10]{Cio} where 2-transitivity is explicitly used). 
We fix a vertex $x \in (\xi_{-}, \xi_{+})$. One has the following Cartan decomposition (see e.g. Ciobotaru \cite[Example 4.10]{Cio}):  $G=K A K$, where 
$K:=G_x= \{g \in G \; \vert \; g(x)=x\}$ and $A:= G_{\xi_-, \xi_{+}}=\{g \in  G \; \vert \; g(\xi_-)=\xi_-, g(\xi_+)=\xi_+\}$.

Note that $A$ is a closed subgroup of $G$ containing $a$. Notice also that each element of $A$ either is elliptic, thus fixing pointwise the axis $(\xi_{-}, \xi_{+})$, or it is hyperbolic and thus translating along the axis $(\xi_{-}, \xi_{+})$. Moreover, every element $g \in G$ is of the form $g=k_1 a^{n}k_2$, for some $k_1,k_2 \in K$ and some $n \in \Z$. For the latter decomposition we used the fact that $a$ has translation length $2$.

By  \cite[Proposition 4.11]{Cio}, we have that 
$$G= \langle G^0_{\xi_+}, G^0_{\xi_-}\rangle.$$

Moreover, by \cite[Proposition 4.15 and Corollary 4.17]{Cio} we have 
$$G_{\xi_{-}}= U_a^{-} A,\;G_{\xi_{+}}= U_a^{+} A,\;G^0_{\xi_-}=U_a^-(A\cap G^0_{\xi_-}) \text{ and } G^0_{\xi_+}=U_a^+(A\cap G^0_{\xi_+}).$$

Notice that $G_{\xi_{-}}$, $G_{\xi_{+}}$, $G^0_{\xi_-}$ and $G^0_{\xi_+}$ are closed subgroups of $G$, and that $A\cap G^0_{\xi_+}=A\cap G^0_{\xi_-}=G^0_{\xi_+}\cap G^0_{\xi_-}$.

The following lemma says that hyperbolic elements in $A$ are boundedly generated by $G_{\xi_{-}}^{0}\cup G_{\xi_{+}}^{0}$.

\begin{Lem}(See the proof of \cite[Lemma 3.5]{CaCi})
\label{lem::bounded_gen}
For every hyperbolic element $\gamma \in A$ there exist $\gamma_1 \in A$ hyperbolic and $u \in (G_{\xi_{+}}^{0} \cap G_{\xi_{-}}^{0})$ such that $\gamma= \gamma_1 u$, $\gamma_1$ has same translation length as $\gamma$ and $\gamma_1$ is a product of $6$ elements from $G_{\xi_{-}}^{0}$ and $G_{\xi_{+}}^{0}$.
\end{Lem}

{\bf Proof:}
As in the proof of \cite[Lemma 3.5]{CaCi} we start with a basic observation. Let $(\eta, \xi_{+})$ be a bi-infinite geodesic line in $T_d$, with $\eta \in \partial T_d \setminus\{\xi_-, \xi_+\}$. Then the intersection $(\eta, \xi_{+})\cap (\xi_{-}, \xi_{+})$ is a geodesic ray $[y, \xi_{+})$, with $y$ a vertex in $(\xi_{-}, \xi_{+})$. We claim that there is some $u \in G_{\xi_+}^{0}$ mapping $(\eta, \xi_{+})$ to $(\xi_{-}, \xi_{+})$ and fixing $[y, \xi_{+})$ pointwise. Indeed, because $\eta$ and $\xi_-$ are opposite to $\xi_{+}$ in $\partial T_d$ and $G$ is $2$-transitive on $\partial T_d$ we obtain that $G_{\xi_+}^{0}$ is transitive on $\partial T_d \setminus \{\xi_{+}\}$ by \cite[Lemma 3.5]{CaCi}; there exists an element $u \in G_{\xi_+}^{0}$ with the desired property. By the same argument applied to the pair of bi-infinite geodesic lines $(\xi_{-}, \eta)$ and $(\xi_{-}, \xi_{+})$, we deduce the existence of $u \in G_{\xi_-}^{0}$ mapping $(\xi_{-}, \eta)$ to $(\xi_{-}, \xi_{+})$ and fixing $(\xi_{-}, \eta) \cap (\xi_{-}, \xi_{+})$ pointwise.

\medskip
Next we claim that for every vertex $y \in (\xi_{-}, \xi_{+})$ there exists $r \in \langle G_{\xi_+}^{0}, G_{\xi_-}^{0} \rangle$, product of three elements from $G_{\xi_{-}}^{0}$ and $G_{\xi_{+}}^{0}$, such that $r$ fixes $x$ and swaps $\xi_-$ to $\xi_+$: i.e. $r(y)=y, r(\xi_-)=\xi_+$ and $r(\xi_+)=\xi_-$. Indeed, fix a vertex $y \in (\xi_{-}, \xi_{+})$. Because $T_d$ is $d$-regular with $d \geq 3$, there exists $\eta \in \partial T_d$ with $(\eta, \xi_+) \cap (\xi_-, \xi_+)= [y, \xi_+)$.  Moreover, we also have that $(\xi_-, \eta) \cap (\xi_-, \xi_+)= [y,\xi_-)$. By the above claim, we can find an element $u \in  G_{\xi_+}^{0} $ fixing $[y, \xi_+)$ pointwise and mapping $[y,\xi_-)$ to $[y,\eta)$. Similarly, there are elements $v, w \in G_{\xi_-}^{0}$ both fixing $[y,\xi_-)$ pointwise and such that $v([y,\eta)) = [y, \xi_+)$ and $w([y, \xi_+)) = u^{-1}([y,\xi_-))$. Now we set $r := vuw$. By construction $r$ fixes the vertex $y$, $r \in \langle G_{\xi_+}^{0}, G_{\xi_-}^{0} \rangle$ and $r$ is the  product of three elements from $G_{\xi_{-}}^{0}$ and $G_{\xi_{+}}^{0}$. Moreover we have
$$r([y, \xi_+)) = vu u^{-1}([y,\xi_-)) = v([y,\xi_-)) = [y,\xi_-)$$
and 
$$r([y,\xi_-)) = vu([y,\xi_-)) = v([y, \eta)) = [y, \xi_+),$$
so that $r$ swaps $[y, \xi_+)$ and $[y,\xi_{-})$. 

Let now $\gamma \in A$ be a hyperbolic element, thus of translation length $2n$, for some $n \in \N$ (recall that $G$ is type-preserving). Fix $x \in (\xi_{-}, \xi_{+})$ and let $y \in (\xi_{-}, \xi_{+})$ be the midpoint of the segment $[x,\gamma(x)]$. Because $\gamma$ has even translation length, $y$ is a vertex of $T_d$. By our second claim above, there exist $r_1, r_2 \in \langle G_{\xi_+}^{0}, G_{\xi_-}^{0} \rangle$, each being product of three elements from $G_{\xi_{-}}^{0}$ and $G_{\xi_{+}}^{0}$, such that $r_1(x)=x$, $r_2(y)=y$ and $r_i(\xi_-)=\xi_+, r_i(\xi_+)=\xi_-$, for $i \in \{1,2\}$. We claim that $r_2 r_1 \in A$ is a hyperbolic element of translation length $2n$ and it is a product of six elements from $G_{\xi_{-}}^{0}$ and $G_{\xi_{+}}^{0}$. Indeed, $r_2 r_1(\xi_-)= r_2(r_1(\xi_-))=r_2(\xi_+)=\xi_-$ and $r_2 r_1(\xi_+)= r_2(r_1(\xi_+))=r_2(\xi_-)=\xi_+$, so $r_2r_1\in A$. Moreover, $r_2 r_1(x)=r_2(x)=\gamma(x)$ and thus $r_2r_1$ is hyperbolic as desired. In particular, we obtain that $u^{-1}:=\gamma^{-1}r_2r_1$ fixes pointwise the bi-infinite geodesic line $(\xi_{-}, \xi_{+})$, thus $u \in G_{\xi_{+}}^{0} \cap G_{\xi_{-}}^{0}$. By taking $\gamma_1:=r_2r_1$ the conclusion follows. 
\hfill$\square$

\begin{Prop}(cf. Cornulier \cite[Proposition 4.1]{Cor})
\label{prop::length_n_proper}
Let $L$ be a length function on $U_a^{+} A$. If $L$ is non-proper on $A$, then $L$ is bounded on $U_a^{+}$ and also on $G^0_{\xi_+}$.
\end{Prop}
{\bf Proof:}
Let $W$ be a compact neighborhood of the identity element $e$ in $G_{\xi_{+}}$, so that $L$ is bounded by a constant $M$ on $W$. 

Suppose that the length function $L$ is not proper on $A$. Then there exists an unbounded sequence $\{a_n\}_{n \in \N} \subset A$ such that $\{L(a_n)\}_{n \in \N}$ is bounded by a constant $M'$.  Then, for every $n \in \N$, we can write $a_n= k_n a^{m_n}k_n'$, for some $k_n,k_n' \in K$ and some $m_n \in \Z$. We obtain  that $\{L(a^{m_n})\}_{n \in \N}$ is bounded by a constant $M''$ and that $m_n \to \infty $ (by extracting a subsequence), when $n \to \infty$. By replacing $m_n$ with $-m_n$ (as $L$ is symmetric) and by taking $u \in U_a^{+}$ one can suppose that $ \lim_{n \to \infty}a^{-m_n}u a^{m_n} =e$. Therefore, for $m_n$ large enough we have that $a^{-m_n}u a^{m_n} \in W$. We obtain that $L(u)\leq M+2M''$, for every $u \in U_a^{+}$. As $G^0_{\xi_+}=U_a^+(A\cap G^0_{\xi_+})$ and because $A\cap G^0_{\xi_+}$ is compact as a closed subgroup of $K$, the function $L$ is also bounded on $G^0_{\xi_+}$. The conclusion follows.
\hfill$\square$

\begin{Cor}
\label{cor::PL}
Let $G$ be a closed non-compact subgroup of $\mathrm{Aut}(T_d)^{+}$ that acts $2$-transitively on $\partial T_d$. Then $G$ has property $PL$.
\end{Cor}
{\bf Proof:}
Let $L$ be a length function on $G$. 
If $L$ is non-proper, then by Proposition \ref{prop::length_n_proper} we have that $L$ is bounded on $G^0_{\xi_+}$ and $G^0_{\xi_-}$. As $G=KAK$ and $L$ is bounded on the compact subgroups $K$, $A\cap G^0_{\xi_+}$ and $A\cap G_{\xi_-}$, it is enough to prove that $L$ remains bounded on the set of all hyperbolic elements in $A$. This follows by applying the bounded generation result of Lemma \ref{lem::bounded_gen} to every hyperbolic element of $A$. Thus $L$ is bounded on $G$ as desired.
\hfill$\square$

\medskip
As we have mentioned above, beside $\mathrm{SL}(2, \Q_p)$, examples of closed non-compact subgroups of $\mathrm{Aut}(T_d)^{+}$ that are topologically simple and act $2$-transitively on $\partial T_d$ are the universal groups $U(F)^{+}$ introduced by Burger--Mozes in \cite[Section~3]{BM00a}. These groups are defined as follows.

\begin{Def}
\label{def::legal_color}
Let $E(T_d)$ be  the set of unoriented edges of the tree $T_d$. Let $\iota: E(T_d) \to \{1,...,d\} $ be a function whose restriction to the star  $E(x)$ of every vertex $x \in T_d$ is   a bijection. Such a function $\iota$  is called a \textbf{legal coloring} of the tree $T_d$.
\end{Def}

\begin{Def}
\label{def::universal_group}
Let $F$ be a subgroup of permutations of the set $\{1,...,d\}$ and let $\iota$ be a legal coloring of $T_d$. The \textbf{universal group}, with respect to $F$ and $\iota$, is defined  as
$$ U(F):= \{g \in \mathrm{Aut}(T_d) \; \vert \; \iota \circ g \circ (\iota \vert_{E(x)})^{-1} \in F, \text{ for every x}\in T_d \}. $$

By $U(F)^{+}$ one denotes the subgroup generated by the edge-stabilizing elements of $U(F)$, and $U(F)^{+} \leq \mathrm{Aut}(T_d)^{+}$. Moreover, Amann \cite[Proposition 52]{Amann} tells us that the group $U(F)$ is independent of the legal coloring $\iota$ of $T_d$.
\end{Def}

Immediately from the definition one deduces that $U(F)$ and $U(F)^{+}$ are closed subgroups of $\mathrm{Aut}(T_d)$. Notice that, when $F$ is the full permutation group $Sym(d)$, then $U(F)= \mathrm{Aut}(T_d)$ and $U(F)^{+}=\mathrm{Aut}(T_d)^{+}$, the latter group being an index $2$, simple subgroup of $\mathrm{Aut}(T_d)$ (for this see Tits \cite{Ti70}).  

An important property of these groups is that $U(F)$ and $U(F)^{+}$ act $2$--transitively on the boundary $\partial T_d$ if and only if $F$ is $2$--transitive. Moreover, $U(F)^{+}$ is either trivial or it is a topologically simple group (see \cite{BM00a, Amann}). 

Moreover, the group $U(F)^{+}$ is not linear. This is because $U(F)^{+}$ has Tits' independence property (see Amann \cite{Amann}); this implies by Caprace-De Medts \cite[Section 2.6]{CaMe11} that the contraction groups $U_{a}^{\pm}$ corresponding to hyperbolic elements $a \in U(F)^{+}$ are not closed. By Wang \cite[Theorem 3.5(ii)]{Wang} we know that the contraction groups corresponding to a $p$-adic Lie group are closed, thus one obtains the non-linearity of $U(F)^{+}$ in characteristic 0. In particular, when $F$ is $2$-transitive, we conclude that the universal group $U(F)^{+}$ is non-linear (in characteristic 0) and has property $PL$.

\noindent
Authors addresses:

\medskip
\noindent
R.T.: Laboratoire de Math\'ematiques d'Orsay\\
Univ. Paris-Sud, CNRS, \\
Universit\'e Paris-Saclay,\\
F-91405 Orsay, France\\
tessera@phare.normalesup.org

\bigskip
\noindent
A.V.: Institut de Math\'ematiques\\
Universit\'e de Neuch\^atel\\
11 Rue Emile Argand\\
CH-2000 Neuch\^atel - Switzerland\\
alain.valette@unine.ch

\bigskip
\noindent
C.C.: Universit\'e de Fribourg\\
Section de math\'ematiques\\
Chemin du Mus\'{e}e 23\\
CH-1700 Fribourg - Switzerland\\
corina.ciobotaru@unifr.ch

\end{document}